\documentclass{amsart}
\usepackage[utf8]{inputenc}

\usepackage[utf8]{inputenc}
\usepackage{amsmath,amssymb,amsfonts,mathtools,amsthm}

\usepackage{mdframed}
\mdfsetup{innerbottommargin=12 pt,innertopmargin=5 pt}

\usepackage{tikz-cd}
\usepackage{tabularray}

\usepackage[hidelinks,pdfencoding=auto,psdextra]{hyperref}
\usepackage{float}

\usepackage{tasks}

\usepackage{array}
\usepackage{siunitx}
\usepackage{booktabs}

\newtheorem{theorem}{Theorem}[section]

\theoremstyle{definition}
\newtheorem{definition}[theorem]{Definition}

\newtheorem{proposition}{Proposition}

\theoremstyle{remark}

\numberwithin{equation}{section}

\newcommand{\Ss}{\mathbb{S}}
\newcommand{\R}{\mathbb{R}}

\newcommand{\So}{\mathrm{SO}}
\newcommand{\SU}{\mathrm{SU}}
\newcommand{\Sl}{\mathrm{SL}(2,\R)}
\newcommand{\Slt}{\mathrm{SL}(3,\R)}
\newcommand{\U}{\mathrm{U}}
\newcommand{\Ima}{\mathrm{Im}}
\newcommand{\iso}{\mathrm{Iso_o}}
\newcommand{\id}{\mathrm{Id}}
\newcommand{\Sp}{\mathrm{Sp}}
\newcommand{\cpt}{\mathbb{C}P^3}
\newcommand{\flag}{F(\mathbb{C}^3)}

\author{Mateo Anarella}
\address{KU Leuven, Department of Mathematics, Celestijnenlaan 200 B – Box 2400, 3001 Leuven, Belgium}
\curraddr{Université Polytechnique Hauts-de-France, Campus Mont Houy 59313 Valenciennes Cedex 9, France.}
\email{mateo.anarella@kuleuven.be}
\thanks{The author is supported by CERAMATHS, University Polytechnique Hauts-de-France / Methusalem grant METH/21/03--long term structural funding of the Flemish Goverment.}
\title{A survey on submanifolds of nearly Kähler spaces}
\begin{document}

\begin{abstract}
In this survey, we give an introduction to nearly Kähler geometry, and list some results on submanifolds of these spaces. This survey tries by no means to be complete.
\end{abstract}

\maketitle
\section{Introduction}
We understand by ``submanifold theory" the study of submanifolds from a Riemannian point of view. Namely, the study of Riemannian immersions, their intrinsic and the extrinsic geometry and the local invariants.
Nonetheless, submanifold theory borrows some concepts from other areas in geometry. For instance, a Lagrangian submanifold is a submanifold for which the symplectic form vanishes everywhere.
Although not every Riemannian manifold is a symplectic manifold, for almost Hermitian manifolds there is a canonical 2-form: the Kähler form.

Gray and Hervella~\cite{grayhervella} divided almost Hermitian manifolds into 16 different classes. Among them we can find Kähler manifolds (Riemannian, symplectic and complex), almost Kähler (Riemannian and symplectic but not complex) and nearly Kähler (only Riemannian). 
This last kind does not have any integrable structure, but it received great interest in the last years. 

% \begin{figure}[ht]
%     \centering
%     \includegraphics[scale=0.4]{diagram.png}
%     \caption{}
%     \label{fig:enter-label}
% \end{figure}

In this survey, I talk about submanifolds in the different nearly Kähler manifolds. Due to works of Nagy and Butruille, mentioned in Section \ref{sectionNK}, we will focus mostly on the four homogeneous ones: $\Ss^6$, $\Ss^3\times\Ss^3$, $\cpt$ and $\flag$.

\section{Nearly Kähler manifolds and splitting theorems}\label{sectionNK}
Let $(N,g,J)$ be an almost Hermitian manifold. That is, $(N,g)$ is a Riemannian manifold and $J$ is an almost complex structure on $N$ compatible with $g$:
\[
J^2=-\id, \ \ \ \ \ \ \ \ \ \ g(JX,JY)=g(X,Y),
\]
for all $X,Y\in\Gamma(TN)$.

\begin{definition}
    A \textbf{nearly Kähler} manifold is an almost Hermitian manifold $(N,g,J)$ such that $\tilde{\nabla} J$ is skew-symmetric. Here, $\tilde{\nabla}$ refers to the Levi-Civita connection associated to $g$. 
    We say that $N$ is a \textbf{strict nearly Kähler manifold} if $\tilde{\nabla}_XJ\neq 0$ for every $X\in\Gamma(TN)$ different from zero.
    
    % If $N$ is not Kähler the we say that $N$ is \textbf{strict nearly Kähler}.
\end{definition}
Although nearly Kähler manifolds are even dimensional, we only find interesting examples from dimension six onwards:
\begin{proposition}
    Two-dimensional and four-dimensional nearly Kähler manifolds are Kähler.
\end{proposition}
Another interesting property of nearly Kähler manifolds is that they are nicely decomposable. Gray and Nagy did significant work in this direction. 
\begin{proposition}[Gray '76 \cite{gray1}, see also \cite{nagyfirst}]
    Let $(N,g,J)$ be a complete and simply connected nearly Kähler manifold. Then $N$ is a Riemannian product $M_1\times M_2$ where $M_1$ is a Kähler manifold and $M_2$ is a strict nearly Kähler manifold.
\end{proposition}

\begin{theorem}[Gray '76 \cite{gray1}]
    There do not exist eight-dimensional strict nearly Kähler manifolds.
\end{theorem}
In fact, this implies that all eight-dimensional nearly Kähler manifolds are either Kähler or a product of a six-dimensional strict nearly Kähler manifold and a two-dimensional Kähler manifold.
In the book of Schäfer~\cite{schafer2} we see that this is true also for pseudo-nearly Kähler manifolds, as long as the three-form $g(\tilde{\nabla}_\cdot J \cdot,\cdot)$ has non-zero length.
\begin{theorem}[Nagy '02 \cite{nagyfirst}]
    A 10-dimensional complete nearly Kähler manifold is either the universal cover of a Riemannian product of a six-dimensional nearly Kähler manifold with a Kähler surface, or the twistor space over a positive, eight-dimensional quaternionic manifold.
\end{theorem}
All eight-dimensional positive quaternionic manifolds are classified, so the 10-dimensional strict nearly Kähler manifolds are as well.

To conclude the discussion on splitting theorems, I present Nagy's theorem, the culmination of these results.
\begin{theorem}[Nagy '02 \cite{nagy}]
    Let $(N,g,J)$ be a complete, simply connected strict nearly Kähler manifold. Then $N$ is a Riemannian product whose factors belong to one of the following three classes:
    \begin{enumerate}
        \item Six-dimensional nearly Kähler manifold,
        \item Homogeneous nearly Kähler manifolds (with some conditions),
        \item Twistor spaces over quaternionic Kähler manifolds with positive scalar curvature.
    \end{enumerate}
\end{theorem}
\noindent The conditions for homogeneous nearly Kähler manifolds in the theorem above are quite involved, for more information see Nagy's original work.

We see in Nagy's theorem that six-dimensional and homogeneous nearly Kähler are among the most important classes.  Butruille classified all nearly Kähler manifolds satisfying both properties:
\begin{theorem}[Butruille '05 \cite{butruille}]
    The only  simply connected homogeneous strict nearly Kähler spaces of dimension six are isometric to $G/H$, where $G$ and $H$ are the Lie groups given in the list.
    \begin{enumerate}
        \item $G=\Ss^3\times\Ss^3\times\Ss^3$, $H=\Delta\Ss^3$. 
        \item $G=G_2$, $H=\SU(3)$. 
        \item $G=\Sp(2)$, $H=\Ss^1\times\SU(2)$. 
        \item $G=\SU(3)$, $H=\Ss^1\times\Ss^1$.
    \end{enumerate}
    The manifolds above are, respectively, $\Ss^3\times\Ss^3$ with a different metric than the product metric, the six-sphere $\Ss^6$ with the canonical metric, the complex projective space $\mathbb{C}P^3$ with a metric different than the Fubini--Study metric, and the manifold of full flags in $\mathbb{C}^3$.
\end{theorem}
Since Tachibana introduced the concept in~\cite{tachibana}, the question if there exist six-dimensional nearly Kähler manifolds which are not homogeneous remained open until Foscolo and Haskins~\cite{foscolo} gave two examples of six-dimensional nearly Kähler manifolds which are not homogeneous, but cohomogeneity-one instead.

Butruille's classification is still open for the pseudo-Riemannian case. Using construction by T-duality, Kath~\cite{kath} gave six analogues of the six-dimensional homogeneous examples:
\begin{tasks}[style=itemize]
% label=(\arabic*)]
(2)
\task $\Sl\times\Sl$, 
\task $\Ss^6_4$, 
\task $\frac{\So^+(2,3)}{\U(1,1)}$, 
\task $\frac{\So^+(4,1)}{\U(2)}$,
\task $\frac{\SU(2,1)}{\U(1)\times \U(1)}$, 
\task $\frac{\Slt}{\R^*\times\So(2)}$.
\end{tasks}
A natural assumption would be that every six-dimensional homogeneous pseud-nearly Kähler manifold is a T-dual of a Riemannian nearly Kähler, but Alekseevsky, Kurglikov and Winther constructed in~\cite{alekseevsky} an example of a pseudo-nearly Kähler six-manifold which is not a T-dual of a Riemannian one.

\section{Submanifold theory}
Let $f:M\to N$ be an isometric immersion of a (pseudo-)Riemannian manifold $M$ into a (pseudo-)Riemannian manifold $N$. 
In submanifold theory, we consider two submanifolds ``equal" if there exists an isometry of the ambient space such that maps one into the other. 
We denote by $\operatorname{Iso}(N)$  the isometry group of $N$, and by $\iso(N)$ the connected component of the identity.

Denote by $\nabla$ and $\tilde{\nabla}$ the Levi-Civita connections of $M$ and $N$, respectively. The connection $\tilde{\nabla}$ splits into a tangent and normal part given by the Gauss and Weingarten formulas below. For simplicity we omit the immersion $f$. For $X,Y\in \Gamma(TM)$ we have that
\[
\tilde{\nabla}_XY=\nabla_XY+h(X,Y), \tag{Gauss formula}
\]
where $h$ is a symmetric $TM^\perp$-valued tensor, called the \textbf{second fundamental form}. For $X\in \Gamma(TM)$ and $\xi\in \Gamma(TM^\perp)$ we have
\[
\tilde{\nabla}_X\xi=-A_\xi X+\nabla^\perp_X\xi, \tag{Weingarten formula}
\]
where $A_\xi$ is symmetric tensor linear in $\xi$, known as \textbf{the shape operator}. The connection $\nabla^\perp\colon \Gamma(TM)\times \Gamma(TM^\perp) \to \Gamma(TM^\perp)$ is known as the \textbf{normal connection}.

We could say that the simplest kind of submanifolds are those with no extrinsic geometry. These are know as totally geodesic submanifolds.
\begin{definition}
    A submanifold is said to be \textbf{totally geodesic} if its geodesics are also geodesics of the ambient space. Equivalently, a submanifold is totally geodesic if and only if the second fundamental form vanishes everywhere.
\end{definition}

An important vector field in submanifold theory that encodes extrinsic curvature is the one known as the mean curvature vector field.
\begin{definition}
    Let $M^m$ be a (pseudo-)Riemannian immersed manifold of a {(pseudo-)}Riemannian manifold $(N^n,g)$, and let $\{e_i\}_1^n$ be a (pseudo-)orthonormal frame with signatures $g(e_i,e_i)=\kappa_i=\pm1$. The \textbf{mean curvature vector field} is defined as     \[H=\frac{1}{n}\sum_{i=1}^n \kappa_i h(e_i,e_i).\]   
\end{definition}
\begin{definition}
    A submanifold is said to be \textbf{minimal} if the mean curvature vector field $H$ vanishes everywhere.
\end{definition}
The name minimal comes from the fact that these submanifolds are critical points of the energy functional.
\begin{definition}
    A submanifold is said to be \textbf{totally umbilical} if for every normal vector field $\xi$ the shape operator $A_\xi$ is a multiple of the identity. Equivalently, there exist a normal vector field $\eta$ such that $h(X,Y)=g(X,Y)\eta$. Moreover, $\eta=H$ in this case.
\end{definition}

As nearly Kähler manifolds are in particular almost Hermitian manifolds, we give an important class of submanifolds involving the complex structure.
\begin{definition}
    Let $(N,g,J)$ be an almost Hermitian manifold. A \textbf{CR-}\textbf{subma\-nifold} $M$ of $N$ is a submanifold such that there exists a holomorphic smooth distribution $D$ on $M$, i.e. $JD\subset D$, and the orthonormal distribution $D^\perp$ is a totally real distribution, i.e. $J(D^\perp)\subset TM^\perp$.

    If $D=TM$, we say that $M$ is an \textbf{almost complex submanifold} (if it is a surface it is also known as pseudo-holomorphic curve or $J$-holomorphic curve). If instead $D$ is the trivial distribution then we say that $M$ is a \textbf{totally real submanifold}. 
\end{definition}
\begin{definition}
    A \textbf{Lagrangian submanifold} $M$ of an almost Hermitian manifold $N$ is a totally real submanifold such that $2\dim M=\dim N$. In the pseudo-Riemannian setting, we add the condition of being non-degenerate.
\end{definition}
\noindent Lagrangian submanifolds of six-dimensional nearly Kähler manifolds are particularly interesting, since they are automatically minimal and orientable.

% \begin{definition}
%     A submanifold $M$ of an almost Hermitian manifold $(N,g,J)$ is said to be \textbf{totally real} if $J$ maps the tangent space into the normal space.
% \end{definition}
% \begin{definition}
%     A totally real submanifold with maximal dimension is called a Lagrangian submanifold. In particular, the dimension of the submanifold is half of the dimension of the ambient space.
% \end{definition}

\begin{definition}
    Given a (pseudo-)Riemannian manifold $N$, a submanifold is \textbf{extrinsically homogeneous} (also known as equivariant) if it is the orbit of a Lie subgroup of $\iso (N)$.
\end{definition}
\begin{definition}
    An action of a Lie group $H$ on a (pseudo-)Riemannian manifold $N$ has \textbf{cohomogeneity one} if $H$ is a Lie subgroup of $\iso(N)$ and the codimension of the principal orbits is one. Note that these principal orbits are extrinsically homogeneous hypersurfaces.
\end{definition}

    Given a hypersurface $M$ in a nearly Kähler manifold $(N,g,J)$ with unit normal~$\xi$, the tangent component of the vector field $JX$ for $X\in\Gamma(TM)$ defines an almost contact structure. Namely, the tensor $\varphi$ given by
    \begin{equation*}
    \varphi (X)=JX-g(JX,\xi)\xi,    
    \end{equation*}
    determines an almost contact structure.

    % From Weingarten formula, we know that $\tilde{\nabla}_X\xi=-AX$. 
    
\begin{definition}
     Let $M$ be a hypersurface in a nearly Kähler manifold $(N,g,J)$, and let $\xi$ be the unit normal vector field. The vector field $J\xi$ is the \textbf{structure vector field} (also known as Reeb vector field).
\end{definition}

\begin{definition}
A Hopf hypersurface is a hypersurface whose structure vector field is a principal vector field, i.e. $A_\xi J\xi=\alpha J\xi$, where $\alpha$ is a smooth function. Equivalently, the integral curves of $J\xi$ are geodesics.
\end{definition}
\begin{definition}
    A hypersurface $M$ of a nearly Kähler manifold $(N,g,J)$ is said to be \textbf{nearly cosymplectic} if $\tilde{\nabla}\varphi$ is skew symmetric.
\end{definition}
\noindent A result by Blair~\cite{blair} says that nearly cosymplectic hypersurfaces satisfy
\begin{equation}
A_\xi\circ\varphi=\varphi\circ A_\xi.  \label{ncos}
\end{equation}
Moreover, all hypersurfaces satisfying \eqref{ncos} are Hopf hypersurfaces.

\section{Submanifolds in the literature}
Because of the splitting theorems by Nagy and Gray, there are not many works on submanifolds of nearly Kähler manifolds of dimension higher than six. However, following a similar construction as for $\Ss^3\times\Ss^3$, we can define a nearly Kähler structure (not necessarilly strict) on $G\times G$, where $G$ is a Lie group.

Aguilar-Suárez and Ruiz-Hernández~\cite{gxg} gave a characterisation of minimal Lagrangian submanifolds of $G\times G$, since in higher dimensions Lagrangian submanifolds are not necessarily minimal. 
Lubbe and Schäfer~\cite{schafer1} studied pseudo-holomorphic curves in $G\times G$. They prove that any totally geodesic pseudo-holomorphic curve is locally determined by a Lie subalgebra $\mathfrak{a}$ of $\mathfrak{g}$, with $0<\dim \mathfrak{a}\leq2$.
\subsection{In six dimensional nearly Kähler manifolds}\label{sectionsixdim}

For general six dimensional nearly Kähler manifolds, Lê and Schwachhöfer~\cite{lorenz} studied the variation of the volume functional for Lagrangian submanifolds, coupled with an analysis of Lagrangian deformations.

Storm~\cite{lagrstorm} constructed examples of Lagrangian submanifolds in $\cpt$ and $\flag$ from minimal submanifolds in their underlying spaces.

In~\cite{loubeau}, Deschamps and Loubeau studied hypersurfaces of both $\mathbb{C}P^3$ and $F(\mathbb{C}^3)$. They proved that there are no hypersurfaces satisfying \eqref{ncos}. Consequently, the only nearly cosymplectic hypersurface of six-dimensional homogeneous nearly Kähler manifolds is $\Ss^5\hookrightarrow\Ss^6$. Moreover, there are no totally geodesic nor totally umbilical hypersurfaces in $\cpt$ or $\flag$.

Hu, Yao and Zhang showed in~\cite{xizejunzeke} that there does not exists a hypersurface in $\Ss^3\times\Ss^3$ for which the shape operator and the almost contact structure tensor~$\varphi$ anticommute. On the other hand, they showed that the only hypersurfaces satisfying such equation in $\Ss^6$, are totally geodesic.

\subsection{In the nearly Kähler \texorpdfstring{$\Ss^6$}{S6}}
The sphere $\Ss^6$ is the most traditional example of a six-dimensional nearly Kähler manifold. It carries the metric inherited from the immersion in $\R^7$ and the almost complex structure is determined by the cross product in the imaginary octonions.

Due to the extensive research in this manifold, to list all the results on submanifolds would be a hard task. Luckily, there are already some surveys that focus particularly on $\Ss^6$. 
For instance, Sekigawa and Hashimoto ~\cite{sekigawasurvey} (2004) wrote about CR-submanifolds, J-holomorphic curves and totally real submanifolds. Dillen also wrote about Lagrangian submanifolds in~\cite{frankisurvey} (1996). Two decades later (2016), Antić and Vrancken also listed works on CR-submanifolds~\cite{anticsurvey}. Thus, we are going to list only works after this date.

In~\cite{madnick} Madnick studied holomorphic curves inside geodesic balls. He proved that if the boundary of the holomorphic curve meets the boundary of the geodesic ball, then the curve is totally geodesic. Along the same lines, Tsouri and Vlachos~\cite{tsourivlachos} studied relations between minimal surfaces and holomorphic curves inside $\Ss^5\hookrightarrow\Ss^6$. Also, Eschenburg and Vlachos~\cite{vlachoseschenburg} gave a characterisation of holomorphic curves in $\Ss^6$ and $\Ss^5$.

Many results regarding Lagrangian submanifolds in $\Ss^6$ have been obtained. Hu, Yin and Yao~\cite{pinching} studied restriction on the Ricci curvature of Lagrangian submanifolds. They conclude that, under some pinching conditions, the only Lagrangian submanifolds are $\Ss^3$ or Berger spheres. Also Lê and Schwachhöfer~\cite{lorenz} studied deformations in Lagrangian submanifolds.

Moruz~\cite{marilena3} studied warped product Lagrangian immersions, finding a relation of such submanifolds with minimal surfaces in $\Ss^5$. Moreover, Ali et al. \cite{alifatemah} studied the homology of warped product Lagrangian submanifolds. In another work by \c{S}ahin \cite{sahin} the topology of submanifolds is also studied.

Enoyoshi and Tsuakada~\cite{emoyoshi} studied Gauss maps and Lagrangian submanifolds: Gauss maps $M\to\tilde{\operatorname{Gr}}_{ass}(\Ima \mathbb{O})$ associated to Lagrangian immersions $M\to\Ss^6$ are harmonic.

Bae, Park and Sekigawa~\cite{baepark} studied quasi contact metric hypersurfaces, and gave a one-parameter family of totally umbilical hypersurfaces.

% J holomorphci~\cite{Bryant}.

% CR submanifods survey~\cite{anticsurvey}.

% On almost complex surfaces of the nearly Kaehler 6-sphere~\cite{opozda2}

% On totally real 3-dimensional submanifolds of the nearly Kaehler 6-sphere~\cite{opozda3}
% Characterization of totally geodesic totally real 3-dimensional submanifolds in the 6-sphere.\cite{anticdjoricvrancken}. 
% Totally real submanifolds chens equality~\cite{dillenvrancken}.
% homogeneous totally real~\cite{mashimo}.
% Totally real surfaces~\cite{opozda}. 

% Classification of totally real submanifolds with K>1/16~\cite{verstraelen}.
% Quasi-Einstein totally real submanifolds of the nearly Kähler 6-sphere\cite{poldillen}.
% Killing vector fields and Lagrangian submanifolds of the nearly Kaehler~\cite{luckilling}.
% The normal curvature of totally real submanifolds of S6(1)~\cite{dillennormal}

% Problems of U simon minimal~\cite{Dillen3}.

\subsection{In the nearly Kähler \texorpdfstring{$\Ss^3\times\Ss^3$}{S3xS3}} This space does not possess the traditional product metric. Instead, the metric on $\Ss^3\times\Ss^3$ is the one that makes the map $\pi:\Ss^3\times\Ss^3\times\Ss^3\to\Ss^3\times\Ss^3$ given by
\[
\pi(a,b,c)=(ac^{-1},bc^{-1}),
\]
into a Riemannian submersion. This space also carries an almost product structure~$P$, which is defined by swapping the tangent vectors at each factor. Most of the authors working on $\Ss^3\times\Ss^3$ make use of this tensor in their classifications.

Djori\'{c}, Djori\'{c} and Moruz~\cite{geodesicsmoruz} gave a parametrization of the geodesics of $\Ss^3\times\Ss^3$ and they provided conditions for when the geodesics (with respect to the nearly Kähler metric) coincide with the geodesics associated to the product metric.

Almost complex surfaces were studied by Bolton et al. in~\cite{dillen}. 
In particular, they classified all totally geodesic almost complex surfaces in $\Ss^3\times\Ss^3$: a flat torus $T^2$ and $\Ss^2(\tfrac{2}{3})$. 
A Bonnet-type theorem for existence an uniqueness is given by Dioos et al. in~\cite{dioos3}.
Moreover, they proved that all flat  almost complex surfaces are part of a family of immersions of $T^2$ into $\Ss^3\times\Ss^3$. 
Dioos in his PhD thesis~\cite{dioosthesis} studied submanifolds of $\Ss^3\times\Ss^3$ in great extent. 
There, he provided a classification of totally geodesic totally real surfaces.
He found that the only one is $\Ss^2(1)$ immersed in one of the factors of $\Ss^3\times\Ss^3$

Anti\'{c} et al.~\cite{marilena1} classified all the surfaces that are invariant under the almost product structure $P$. 
They found that they are either the torus $T^2$ of the classifications in~\cite{dillen} and~\cite{dioos3}, any  surface in $p:N^2\to\Ss^3$ immersed by $p\mapsto(p,p)$, or solutions of a specific differential equation. 
Furthermore, they classified all slant $P$-invariant surfaces (surfaces with constant Kähler angle). 
This last condition facilitates the solution of the aforementioned differential equation.

Dioos et al.~\cite{wang} proved that any Lagrangian submanifold of a six-dimensional nearly Kähler manifold with parallel second fundamental form is
totally geodesic (this is also true for the pseudo-Riemannian case). 
Moreover, they provided a full classification of totally geodesic Lagrangian submanifolds of $\Ss^3\times\Ss^3$. 
Although they listed six examples, because of~\cite{propertiess3s3} these are only two.
They are either a round sphere $\Ss^3(\tfrac{3}{4})$ or a Berger sphere.

In~\cite{dioos2}, Dioos, Vrancken and Wang classified all Lagrangian submanifolds with constant sectional curvature: the already mentioned totally geodesic sphere $\Ss^3(\tfrac{3}{4})$, a flat torus $T^3$ and $\R P^3(\tfrac{3}{16})$.

Bekta\c{s} et al.~\cite{bektas} studied Lagrangian submanifolds where the angles between $PE_i$ and $E_i$ are constant, for a specific frame $\{E_i\}_i$. They found that they must be either a totally geodesic submanifolds or a submanifold with constant sectional curvature. 
As a consequence, the authors classified  all extrinsically homogeneous Lagrangian submanifolds, since these have constant angles functions.

% Even though is not specified in the article, they classified all extrinsically homogeneous Lagrangian submanifolds. Using an argument of unicity on the frame, we can see that for extrinsically homogeneous Lagrangian submanifolds these angles are constant.

We can find examples of non-extrinsically homogeneous Lagrangian submanifolds in~\cite{moruz}. There, the same authors as in~\cite{bektas} showed a way to construct Lagrangian submanifolds from a minimal submanifold immersed in $\Ss^3$.

Three-dimensional CR-submanifolds have been studied by Anti\'{c} et al in~\cite{moruz2}.
They considered those CR-submanifolds whose associated distributions are totally geodesic.
As a continuation, the authors studied in~\cite{natasja} how the almost product structure $P$ acts on the distributions of three dimensional CR-submanifolds.

Hu, Yao and Zhang~\cite{firsthyperzeke} showed that there are no totally umbilical hypersurfaces nor hypersurfaces with parallel second fundamental form in $\Ss^3\times\Ss^3$. Furthermore, they classified all hypersurfaces that satisfy \eqref{ncos} and have three different principal curvatures. In~\cite{zejunzeke1}, Hu and Yao proved that there does not exist a Hopf hypersurface with only two different principal curvatures. This motivated the authors~\cite{zejunzeke1,zejunzeke2} to study Hopf hypersurfaces with at least three different principal curvatures. Additionally< they require that $P$ preserves the orthogonal space to $J\xi$, where $\xi$ is the unit normal. They found a family of immersions from $\Ss^3\times\Ss^2$ into $\Ss^3\times\Ss^3$ which are Hopf hypersurfaces.

Hu et al.~\cite{Zeke} proved that there are no locally conformally flat hypersurfaces nor Einstein Hopf hypersurfaces. Moreover, they provided a Simons type integral inequality for minimal hypersurfaces.

In the same line as the previous articles, Chen and Yang~\cite{chenyang} studied hypersurfaces satisfying a similar equation to \eqref{ncos}, replacing the shape operator by the Ricci tensor. Adding the condition of being Hopf, they proved that the mean curvature is $2\alpha$, where $\alpha$ is the principal curvature associated to the structure vector field. They also show that the mean curvature is constant if $P$ preserves the normal space to the structure vector field.
Hu, Yao and Zhang~\cite{xizejun} also proved that the following conditions are equivalent:
\begin{itemize}
    \item $P$ preserving the normal spaces to the structure vector field,
    \item the associated principal curvature being constant,
    \item the hypersurface being congruent to a product $\Ss^3\times\Sigma$
\end{itemize}
 with $\Sigma\hookrightarrow\Ss^3$. This generalised a previous result made by Djori\'c~\cite{milos}.

\subsection{In the nearly Kähler \texorpdfstring{$\mathbb{C}P^3$}{CP3}} The nearly Kähler complex projective space does not carry the traditional Fubini--Study metric. Instead, the nearly Kähler metric arises from the fact that $\cpt$ is the twistor space over the quaternionic Kähler $\mathbb{H}P^1\cong \Ss^4$. Namely, $\cpt$ is a fibre bundle by 2-spheres over $\Ss^4$. Then, there exist two distinguished distributions $\mathcal{D}_1^2$ and $\mathcal{D}_2^4$. By modifying complex structure on $\mathcal{D}_1$ and the Kähler metric on $\mathcal{D}_2$, we get a nearly Kähler manifold.

It is clear that the flag manifold and $\mathbb{C}P^3$ are the least studied of the homogeneous six-dimensional nearly Kähler spaces. Xu studied $J$-holomorphic curves~\cite{xufeng}. Also, Aslan studied them in~\cite{aslan2}.

Aslan~\cite{Aslan} and Liefsoens~\cite{liefsoens} studied separately extrinsically homogeneous and totally geodesic Lagrangian submanifolds. Among them, we find $\Ss^1\times\Ss^2$, a Berger sphere, and $\R P^3$, the last one being the only totally geodesic Lagrangian submanifold.

For hypersurfaces see the work by Deschamps and Loubeau in Subsection \ref{sectionsixdim}.

\subsection{In the nearly Kähler \texorpdfstring{$F(\mathbb{C}^3)$}{F(C3)}}
The manifold of full flags in $\mathbb{C}^3$ carries a nearly Kähler structure similar to $\cpt$. It is also a twistor space, but over $\mathbb{C}P^2$. So in the same way, from a Kähler structure we construct a nearly Kähler structure.

Vrancken and Cwiklinski classified  in~\cite{kwilinski} all the totally geodesic almost complex surfaces. There, they found that there are 3 such surfaces, all spaces forms: a sphere of radius 1, a sphere of radius $\tfrac{1}{2}$ and the flat torus. 

Storm classified in~\cite{reinierstorm} all totally geodesic Lagrangian submanifolds. These are either the manifold of full flags in $\R^3$ or a Berger sphere. Moreover, he classified all extrinsically homogeneous Lagrangian submanifolds. He showed that in addition to the previous examples we have $\R P^3$, with a Berger-like metric.

For hypersurfaces see the work by Deschamps and Loubeau in Subsection \ref{sectionsixdim}.

\subsection{In pseudo-nearly Kähler manifolds}
The study of submanifolds in six-dimensional homogeneous pseudo-nearly Kähler manifolds has just started. In his book, Schäfer~\cite{schafer2} did great work on general theory of Lagrangian submanifolds.

Similar to the classification in $\Ss^3\times\Ss^3$ in~\cite{dillen}, Ghandour and Vrancken~\cite{ghandour} classified almost complex totally geodesic surfaces in $\Sl\times\Sl$. These are a flat torus $T^2$, $\R^2$ or $H^2(-\tfrac{4}{3})$. Anarella and Dekimpe~\cite{dekimpe} proved instead that the only two totally real totally geodesic surfaces are $H^2_1(-\tfrac{3}{2})$ and $H^2(-\tfrac{3}{2})$. Dekimpe also studied degenerate almost complex surfaces in~\cite{dekimpe2}. 

When it comes to Lagrangian submanifolds of $\Sl\times\Sl$, Anarella and Van der Veken~\cite{anarella} found a full classification up to congruence. There are three non-congruent totally geodesic Lagrangian submanifolds, $H^3_1(-\tfrac{3}{2})$, and two anti-de Sitter spaces with different Berger-like metrics.

In the same paper~\cite{kwilinski} as for $\flag$, Cwiklinski and Vrancken classified all totally geodesic surfaces of $\tfrac{\SU(2,1)}{\text{U} (1)\times \text{U}(1)}$. Analogous to the Riemannian case they found three surfaces: a sphere $\Ss^2(\tfrac{1}{2})$, $H^2(-4)$ and $H^2(-1)$.

\section{Personal comments}
To finish, I would like to point out the similarity of the works for $\Ss^3\times\Ss^3$~\cite{bektas}, $\cpt$ and $\flag$~\cite{lagrstorm}, and  $\Ss^6$~\cite{marilena3}. In all of these works, the authors construct Lagrangian submanifolds from minimal surfaces immersed in a smaller space. These are the only examples of non-extrinsically homogeneous Lagrangian submanifolds I am aware of. 

It would be interesting to know if we can construct this kind of Lagrangian submanifolds for any nearly Kähler space, if they are always non-extrinsically homogeneous, and if they are the only non-extrinsically homogeneous Lagrangian submanifolds.

\bibliographystyle{abbrv}
 % \bibliographystyle{alpha}.
\footnotesize
\bibliography{survey}

\begin{thebibliography}{10}

\bibitem{gxg}
R.~Aguilar-Su\'{a}rez and G.~Ruiz-Hern\'{a}ndez.
\newblock A characterization of minimal {L}agrangian submanifolds of the nearly {K}\"{a}hler {$G \times G$}.
\newblock {\em Bull. Belg. Math. Soc. Simon Stevin}, 30(2):140--162, 2023.

\bibitem{alekseevsky}
D.~V. Alekseevsky, B.~S. Kruglikov, and H.~Winther.
\newblock Homogeneous almost complex structures in dimension 6 with semi-simple isotropy.
\newblock {\em Ann. Global Anal. Geom.}, 46(4):361--387, 2014.

\bibitem{alifatemah}
A.~Ali, F.~Mofarreh, N.~Alluhaibi, and P.~Laurian-Ioan.
\newblock Null homology in warped product {L}agrangian submanifolds of the nearly {K}aehler {$\mathbb{S}^6$} and its applications.
\newblock {\em J. Geom. Phys.}, 158:103859, 13, 2020.

\bibitem{anarella}
M.~Anarella and J.~V. der Veken.
\newblock Totally geodesic lagrangian submanifolds of the pseudo-nearly k\"ahler $\mathrm{SL}(2,\mathbb{R})\times\mathrm{SL}(2,\mathbb{R})$.
\newblock {\em arXiv preprint 2307.13389}, 2023.

\bibitem{natasja}
M.~Anti\'{c}, N.~Djurdjevi\'{c}, and M.~Moruz.
\newblock C{R} submanifolds of the nearly {K}\"{a}hler {$\mathbb S^3\times\mathbb S^3$} characterised by properties of the almost product structure.
\newblock {\em Mediterr. J. Math.}, 15(3):Paper No. 111, 28, 2018.

\bibitem{moruz2}
M.~Anti\'{c}, N.~Djurdjevi\'{c}, M.~Moruz, and L.~Vrancken.
\newblock Three-dimensional {CR} submanifolds of the nearly {K}\"{a}hler {$\mathbb{S}^3\times \mathbb{S}^3$}.
\newblock {\em Ann. Mat. Pura Appl. (4)}, 198(1):227--242, 2019.

\bibitem{marilena1}
M.~Anti\'{c}, Z.~Hu, M.~Moruz, and L.~Vrancken.
\newblock Surfaces of the nearly {K}\"{a}hler {$\mathbb S^3\times\mathbb S^3$} preserved by the almost product structure.
\newblock {\em Math. Nachr.}, 294(12):2286--2301, 2021.

\bibitem{anticsurvey}
M.~Anti\'{c} and L.~Vrancken.
\newblock C{R}-submanifolds of the nearly {K}\"{a}hler 6-sphere.
\newblock In {\em Geometry of {C}auchy-{R}iemann submanifolds}, pages 57--90. Springer, Singapore, 2016.

\bibitem{aslan2}
B.~Aslan.
\newblock Transverse {$J$}-holomorphic curves in nearly {K}\"{a}hler {$\mathbb{CP}^3$}.
\newblock {\em Ann. Global Anal. Geom.}, 61(1):115--157, 2022.

\bibitem{Aslan}
B.~Aslan.
\newblock Special {L}agrangians in nearly {K}\"{a}hler {$\mathbb{CP}^3$}.
\newblock {\em J. Geom. Phys.}, 184:Paper No. 104713, 18, 2023.

\bibitem{baepark}
J.~Bae, J.~Park, and K.~Sekigawa.
\newblock A one-parameter family of totally umbilical hyperspheres in the nearly {K}\"{a}hler 6-sphere.
\newblock {\em J. Korean Math. Soc.}, 55(4):963--974, 2018.

\bibitem{bektas}
B.~Bekta\c{s}, M.~Moruz, J.~Van~der Veken, and L.~Vrancken.
\newblock Lagrangian submanifolds with constant angle functions of the nearly {K}\"{a}hler {$\mathbb{S}^3\times\mathbb{S}^3$}.
\newblock {\em J. Geom. Phys.}, 127:1--13, 2018.

\bibitem{moruz}
B.~Bekta\c{s}, M.~Moruz, J.~Van~der Veken, and L.~Vrancken.
\newblock Lagrangian submanifolds of the nearly {K}\"{a}hler {$\mathbb S^3 \times \mathbb S^3$} from minimal surfaces in {$\mathbb S^3$}.
\newblock {\em Proc. Roy. Soc. Edinburgh Sect. A}, 149(3):655--689, 2019.

\bibitem{blair}
D.~E. Blair.
\newblock Almost contact manifolds with {K}illing structure tensors.
\newblock {\em Pacific J. Math.}, 39:285--292, 1971.

\bibitem{dillen}
J.~Bolton, F.~Dillen, B.~Dioos, and L.~Vrancken.
\newblock Almost complex surfaces in the nearly {K}\"{a}hler {$\Ss^3\times\Ss^3$}.
\newblock {\em Tohoku Math. J. (2)}, 67(1):1--17, 2015.

\bibitem{butruille}
J.-B. Butruille.
\newblock Classification des vari\'{e}t\'{e}s approximativement k\"{a}hleriennes homog\`enes.
\newblock {\em Ann. Global Anal. Geom.}, 27(3):201--225, 2005.

\bibitem{chenyang}
X.~Chen and Y.~Yang.
\newblock Hopf hypersurfaces of the homogeneous nearly {K}\"{a}hler {$\mathbb{S}^3\times\mathbb{S}^3$} satisfying certain commuting conditions.
\newblock {\em Bull. Korean Math. Soc.}, 59(6):1567--1594, 2022.

\bibitem{sahin}
F.~\c{S}ahin.
\newblock Homology of submanifolds of six dimensional sphere.
\newblock {\em J. Geom. Phys.}, 145:103471, 6, 2019.

\bibitem{kwilinski}
K.~Cwiklinski and L.~Vrancken.
\newblock Almost complex surfaces in the nearly {K}\"{a}hler flag manifold.
\newblock {\em Results Math.}, 77(3):Paper No. 134, 17, 2022.

\bibitem{dekimpe}
K.~Dekimpe.
\newblock {\em Contributions to the theory of surfaces in almost complex spaces}.
\newblock PhD Thesis. KU Leuven, 2023.

\bibitem{dekimpe2}
K.~Dekimpe.
\newblock Degenerate almost complex surfaces in the nearly k\"ahler $\mathrm{SL}_2\mathbb{R}\times \mathrm{SL}_2\mathbb{R}$.
\newblock {\em arXiv preprint}, 2023.

\bibitem{loubeau}
G.~Deschamps and E.~Loubeau.
\newblock Hypersurfaces of the nearly {K}\"{a}hler twistor spaces {$\mathbb C{\rm P}^3$} and {$\mathbb F_{1,2}$}.
\newblock {\em Tohoku Math. J. (2)}, 73(4):627--642, 2021.

\bibitem{frankisurvey}
F.~Dillen.
\newblock Lagrangian submanifolds of {$\Ss^6(1)$}.
\newblock In {\em Geometry and topology of submanifolds, {VIII} ({B}russels, 1995/{N}ordfjordeid, 1995)}, pages 138--144. World Sci. Publ., River Edge, NJ, 1996.

\bibitem{dioosthesis}
B.~Dioos.
\newblock {\em Submanifolds of the nearly {K}ähler manifold {$\Ss^3\times\Ss^3$}}.
\newblock PhD thesis. KU Leuven, 2015.

\bibitem{dioos3}
B.~Dioos, H.~Li, H.~Ma, and L.~Vrancken.
\newblock Flat almost complex surfaces in the homogeneous nearly {K}\"{a}hler {$\Ss^3\times \Ss^3$}.
\newblock {\em Results Math.}, 73(1):Paper No. 38, 24, 2018.

\bibitem{dioos2}
B.~Dioos, L.~Vrancken, and X.~Wang.
\newblock Lagrangian submanifolds in the homogeneous nearly {K}\"{a}hler {$\mathbb S^3\times\mathbb S^3$}.
\newblock {\em Ann. Global Anal. Geom.}, 53(1):39--66, 2018.

\bibitem{milos}
M.~Djori\'{c}.
\newblock Hypersurfaces of the homogeneous nearly {K}\"{a}hler {$\mathbb S^3 \times \mathbb S^3$} whose normal vector field is {$\mathcal P$}-principal.
\newblock {\em Mediterr. J. Math.}, 18(6):Paper No. 251, 20, 2021.

\bibitem{geodesicsmoruz}
M.~Djori\'{c}, M.~Djori\'{c}, and M.~Moruz.
\newblock Geodesic lines on nearly {K}\"{a}hler {$\mathbb{S}^3\times\mathbb{S}^3$}.
\newblock {\em J. Math. Anal. Appl.}, 466(1):1099--1108, 2018.

\bibitem{emoyoshi}
K.~Enoyoshi and K.~Tsukada.
\newblock Lagrangian submanifolds of {$\Ss^6$} and the associative {G}rassmann manifold.
\newblock {\em Kodai Math. J.}, 43(1):170--192, 2020.

\bibitem{vlachoseschenburg}
J.-H. Eschenburg and T.~Vlachos.
\newblock Pseudoholomorphic curves in {$\mathbb S^6$} and {$\mathbb S^5$}.
\newblock {\em Rev. Un. Mat. Argentina}, 60(2):517--537, 2019.

\bibitem{foscolo}
L.~Foscolo and M.~Haskins.
\newblock New {$G_2$}-holonomy cones and exotic nearly {K}\"{a}hler structures on {$\Ss^6$} and {$\Ss^3\times \Ss^3$}.
\newblock {\em Ann. of Math. (2)}, 185(1):59--130, 2017.

\bibitem{ghandour}
E.~Ghandour and L.~Vrancken.
\newblock Almost complex surfaces in the nearly kähler {$\Sl\times\Sl$}.
\newblock {\em Mathematics}, 8(7), 2020.

\bibitem{gray1}
A.~Gray.
\newblock The structure of nearly {K}\"{a}hler manifolds.
\newblock {\em Math. Ann.}, 223(3):233--248, 1976.

\bibitem{grayhervella}
A.~Gray and L.~M. Hervella.
\newblock The sixteen classes of almost {H}ermitian manifolds and their linear invariants.
\newblock {\em Ann. Mat. Pura Appl. (4)}, 123:35--58, 1980.

\bibitem{sekigawasurvey}
H.~Hashimoto and K.~Sekigawa.
\newblock Submanifolds of a nearly {K}\"{a}hler 6-dimensional sphere.
\newblock In {\em Proceedings of the {E}ighth {I}nternational {W}orkshop on {D}ifferential {G}eometry}, pages 23--45. Kyungpook Nat. Univ., Taegu, 2004.

\bibitem{Zeke}
Z.~Hu, M.~Moruz, L.~Vrancken, and Z.~Yao.
\newblock On the nonexistence and rigidity for hypersurfaces of the homogeneous nearly {K}\"{a}hler {$\mathbb S^3\times \mathbb S^3$}.
\newblock {\em Differential Geom. Appl.}, 75:Paper No. 101717, 22, 2021.

\bibitem{zejunzeke1}
Z.~Hu and Z.~Yao.
\newblock On {H}opf hypersurfaces of the homogeneous nearly {K}\"{a}hler {$\Ss^3\times \Ss^3$}.
\newblock {\em Ann. Mat. Pura Appl. (4)}, 199(3):1147--1170, 2020.

\bibitem{pinching}
Z.~Hu, Z.~Yao, and J.~Yin.
\newblock On {R}icci curvature pinching of {L}agrangian submanifolds in the homogeneous nearly {K}\"{a}hler {$\mathbb S^6(1)$}.
\newblock {\em Results Math.}, 75(2):Paper No. 52, 7, 2020.

\bibitem{xizejunzeke}
Z.~Hu, Z.~Yao, and X.~Zhang.
\newblock Hypersurfaces of the homogeneous nearly {K}\"{a}hler {$\mathbb{S}^6$} and {$\mathbb{S}^3\times\mathbb{S}^3$} with anticommutative structure tensors.
\newblock {\em Bull. Belg. Math. Soc. Simon Stevin}, 26(4):535--549, 2019.

\bibitem{firsthyperzeke}
Z.~Hu, Z.~Yao, and Y.~Zhang.
\newblock On some hypersurfaces of the homogeneous nearly {K}\"{a}hler {$\Ss^3\times\Ss^3$}.
\newblock {\em Math. Nachr.}, 291(2-3):343--373, 2018.

\bibitem{kath}
I.~Kath.
\newblock Killing spinors on pseudo-{Riemannian} manifolds.
\newblock {\em Habilitationsschrift an der Humboldt-Universität zu Berlin}, 1999.

\bibitem{lorenz}
H.~V. L\^{e} and L.~Schwachh\"{o}fer.
\newblock Lagrangian submanifolds in strict nearly {K}\"{a}hler 6-manifolds.
\newblock {\em Osaka J. Math.}, 56(3):601--629, 2019.

\bibitem{liefsoens}
M.~Liefsoens.
\newblock {\em Nearly {K}ähler {$\cpt$}: From the {H}opf Fibration to {L}agrangian Submanifolds}.
\newblock Msc thesis. KU Leuven. Faculteit Wetenschappen, Leuven, 2022.

\bibitem{schafer1}
F.~Lubbe and L.~Sch\"{a}fer.
\newblock Pseudo-holomorphic curves in nearly {K}\"{a}hler manifolds.
\newblock {\em Differential Geom. Appl.}, 36:24--43, 2014.

\bibitem{madnick}
J.~Madnick.
\newblock Free-boundary problems for holomorphic curves in the 6-sphere.
\newblock {\em Math. Z.}, 303(4):Paper No. 85, 17, 2023.

\bibitem{marilena3}
M.~Moruz.
\newblock Lagrangian warped product immersions in {$\mathbb S^6$}.
\newblock {\em Arch. Math. (Basel)}, 113(3):325--336, 2019.

\bibitem{propertiess3s3}
M.~Moruz and L.~Vrancken.
\newblock Properties of the nearly {K}\"{a}hler {$\mathbb S^3\times\mathbb S^3$}.
\newblock {\em Publ. Inst. Math. (Beograd) (N.S.)}, 103(117):147--158, 2018.

\bibitem{nagy}
P.-A. Nagy.
\newblock Nearly {K}\"{a}hler geometry and {R}iemannian foliations.
\newblock {\em Asian J. Math.}, 6(3):481--504, 2002.

\bibitem{nagyfirst}
P.-A. Nagy.
\newblock On nearly-{K}\"{a}hler geometry.
\newblock {\em Ann. Global Anal. Geom.}, 22(2):167--178, 2002.

\bibitem{schafer2}
L.~Sch\"{a}fer.
\newblock {\em Nearly pseudo-{K}\"{a}hler manifolds and related special holonomies}, volume 2201 of {\em Lecture Notes in Mathematics}.
\newblock Springer, Cham, 2017.

\bibitem{reinierstorm}
R.~Storm.
\newblock Lagrangian submanifolds of the nearly {K}\"{a}hler full flag manifold {$F_{1,2}(\mathbb C^3)$}.
\newblock {\em J. Geom. Phys.}, 158:103844, 16, 2020.

\bibitem{lagrstorm}
R.~Storm.
\newblock A note on {L}agrangian submanifolds of twistor spaces and their relation to superminimal surfaces.
\newblock {\em Differential Geom. Appl.}, 73:101669, 10, 2020.

\bibitem{tachibana}
S.-i. Tachibana.
\newblock On almost-analytic vectors in almost-{K}\"{a}hlerian manifolds.
\newblock {\em Tohoku Math. J. (2)}, 11:247--265, 1959.

\bibitem{tsourivlachos}
A.-S. Tsouri and T.~Vlachos.
\newblock Minimal surfaces in spheres and a {R}icci-like condition.
\newblock {\em Manuscripta Math.}, 166(3-4):561--588, 2021.

\bibitem{xufeng}
F.~Xu.
\newblock Pseudo-holomorphic curves in nearly {K}\"{a}hler {$\cpt$}.
\newblock {\em Differential Geom. Appl.}, 28(1):107--120, 2010.

\bibitem{zejunzeke2}
Z.~Yao and Z.~Hu.
\newblock On {H}opf hypersurfaces of the homogeneous nearly {K}\"{a}hler {$\Ss^3\times\Ss^3$}, {II}.
\newblock {\em Manuscripta Math.}, 168(3-4):371--402, 2022.

\bibitem{xizejun}
Z.~Yao, X.~Zhang, and Z.~Hu.
\newblock Hypersurfaces of the homogeneous nearly {K}\"{a}hler {$\Ss^3\times \Ss^3$} with {$P$}-invariant holomorphic distributions.
\newblock {\em J. Geom. Anal.}, 32(7):Paper No. 209, 16, 2022.

\bibitem{wang}
Y.~Zhang, B.~Dioos, Z.~Hu, L.~Vrancken, and X.~Wang.
\newblock Lagrangian submanifolds in the 6-dimensional nearly {K}\"{a}hler manifolds with parallel second fundamental form.
\newblock {\em J. Geom. Phys.}, 108:21--37, 2016.

\end{thebibliography}

\end{document}